\documentclass[dvips,a4paper]{amsproc}
\usepackage{amsmath, amssymb, amsfonts, amscd}

\usepackage{graphicx, color}
\usepackage[hypertex]{hyperref}
\usepackage{tikz}
\usepackage[all]{xy}
\numberwithin{equation}{section}
\theoremstyle{plain}
 \newtheorem{thm}{Theorem}[section]
 \newtheorem{cor}[thm]{Corollary}
 \newtheorem{lem}[thm]{Lemma}
\theoremstyle{definition}
 \newtheorem{defn}[thm]{Definition}
 \newtheorem{exmp}[thm]{Example}
\theoremstyle{remark}
 \newtheorem{rem}[thm]{Remark}
\DeclareMathOperator{\trace}{trace}
\DeclareMathOperator{\supp}{supp}
\DeclareMathOperator{\codim}{codim}
\DeclareMathOperator{\ord}{ord}
\DeclareMathOperator{\idx}{idx}

\def\D{\texttt{D}}
\def\Idx{\texttt{Idx}}
\def\Ix{\texttt{Ix}}
\def\IxF{\texttt{IxF}}
\def\L{\texttt{L}}
\def\IL{\texttt{IL}}
\begin{document}
\title[star-shaped Kac-Moody root systems]
{Algorithm classifying roots of \\ star-shaped Kac-Moody root systems}
\author{Toshio \textsc{Oshima}}
\email{oshima@ms.u-tokyo.ac.jp}

\begin{abstract}
There is a correspondence between the imaginary roots of a star-shaped Kac-Moody root system 
and the spectral types of non-rigid Fuchsian ordinary differential equations. 
The norm of a root corresponds to the index of rigidity of the equation, and the action of the Weyl group corresponds to the middle convolution of the equation. 
It is known that the Weyl group has finite orbits in the set of imaginary roots with a given norm.
The purpose of this paper is to present an algorithm implemented in computer algebra that 
effectively computes representatives of these orbits. 
These representatives can be applied to the classification and construction of higher-dimensional 
Painlev\'e-type equations. 
We explain the relationship between the roots of a star-shaped Kac-Moody root system and these equations, along with the related computer programs.
\end{abstract}
\maketitle
\section{Introduction}
We first review a star-shaped Kac-Moody root system.
Let
\begin{equation}
  I:=\{0,\,(j,\nu)\,;\,j=0,1,\ldots,\ \nu=1,2,\ldots\}.
\end{equation}
be a set of indices and 
let $\mathfrak h$ be an infinite dimensional real vector space with 
the set of basis $\Pi$, where
\begin{equation}
 \Pi=\{\alpha_i\,;\,i\in I\}
    =\{\alpha_0,\ \alpha_{j,\nu}\,;\,j=0,1,2,\ldots,\ \nu=1,2,\ldots\}
{.}
\end{equation}
Put
\begin{align}
 I' &:= I\setminus\{0\},\qquad \Pi':=\Pi\setminus\{\alpha_0\},\\
 Q  &:=\sum_{\alpha\in\Pi}\mathbb Z\alpha\ \supset\ 
 Q_+:=\sum_{\alpha\in\Pi}\mathbb Z_{\ge0}\alpha,\ \ Q_-:=-Q_+.
\end{align}
Then $\alpha\in Q$ has the expression
\begin{equation}\label{eq:root}
  \alpha=n\alpha_0
   + \sum_{j=0}^{p-1}\sum_{\nu=1}^{n_j-1} n_{j,\nu}\alpha_{j,\nu}
\end{equation}
with $n_{j,\nu}\in\mathbb Z$, and $\alpha\in Q_+$ if and only if $n_{j,\nu}\in\mathbb Z_{\ge 0}$.
Moreover $\alpha\in Q$ is {\sl indivisible} means that $\frac1k\alpha\not\in Q$ for $k=2,3,\ldots$.  
We define an indefinite symmetric bilinear form on $\mathfrak h$ by\\[-15pt]
\begin{equation}\label{eq:PIKac}%
 \begin{split}
 (\alpha|\alpha)
 &= 2\qquad\ \,(\alpha\in\Pi),\\
 (\alpha_0|\alpha_{j,\nu})
 &=-\delta_{\nu,1},\\
 (\alpha_{i,\mu}|\alpha_{j,\nu})
 &=\begin{cases}
    0 &(i\ne j\text{ \ \ or \ \ }|\mu-\nu|>1),\\
    -1&(i=j\text{ \ and \ }|\mu-\nu|=1).
  \end{cases}
 \end{split}\quad\ \ \raisebox{10pt}{\text{\small
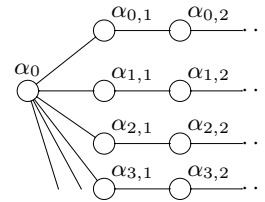
\begin{xy}
\ar@{-}               *++!D{\text{$\alpha_0$}}  *\cir<4pt>{}="O";
             (10,0)   *+!L!D{\text{$\alpha_{1,1}$}} *\cir<4pt>{}="A",
\ar@{-} "A"; (20,0)   *+!L!D{\text{$\alpha_{1,2}$}} *\cir<4pt>{}="B",
\ar@{-} "B"; (30,0)   *{\cdots}, 
\ar@{-} "O"; (10,-7)  *+!L!D{\text{$\alpha_{2,1}$}} *\cir<4pt>{}="C",
\ar@{-} "C"; (20,-7)  *+!L!D{\text{$\alpha_{2,2}$}} *\cir<4pt>{}="E",
\ar@{-} "E"; (30,-7)  *{\cdots}
\ar@{-} "O"; (10,8)   *+!L!D{\text{$\alpha_{0,1}$}} *\cir<4pt>{}="D",
\ar@{-} "D"; (20,8)   *+!L!D{\text{$\alpha_{0,2}$}} *\cir<4pt>{}="F",
\ar@{-} "F"; (30,8)   *{\cdots}
\ar@{-} "O"; (10,-13) *+!L!D{\text{$\alpha_{3,1}$}} *\cir<4pt>{}="G",
\ar@{-} "G"; (20,-13) *+!L!D{\text{$\alpha_{3,2}$}} *\cir<4pt>{}="H",
\ar@{-} "H"; (30,-13) *{\cdots},
\ar@{-} "O"; (7, -13),
\ar@{-} "O"; (4, -13),
\end{xy}}}
\end{equation}
We represent this relation by a diagram which 
contains small circles representing elements of $\Pi$ and 
connected two circles representing $\alpha,\,\beta\in \Pi$ by a line if $(\alpha|\beta)=-1$. 

The star-shaped root system is the set $\Pi$ with the inner products of elements of $\Pi$
defined by the above $\langle\ ,\ \rangle$. 
An elements of $\Pi$ is called a \textsl{simple roots} of the root system 
and the \textsl{Weyl group} $W$ of the root system is 
generated by the 
\textsl{simple reflections}
\begin{equation}\label{eq:sref}
 s_\alpha\,:\,\mathfrak h\ni x\mapsto
 x-(x|\alpha)\alpha\in\mathfrak h
\end{equation}
with $\alpha \in \Pi$. 
The group generated by linear transformations $\sigma_{i,j}$ of $\mathfrak h$
with 
\[\sigma_{i,j}(\alpha_{i,\nu})=\alpha_{j,\nu},\ \sigma_{i,j}(\alpha_{j,\nu})=\alpha_{i,\nu}, \
\sigma_{i,j}(\alpha_{k,\nu})=\alpha_{k,\nu}\ \ (k\ne i,j),\  \sigma_{i,j}(\alpha_0)=\alpha_0
\]
for $0\le i<j$ is denoted by $S_\infty$.
The group generated by $W$ and $S_\infty$ is denote by $\widetilde W$, which
keeps $(\ |\ )$ invariant.
\begin{rem}[\cite{Kc}]\label{rem:Kac}
The set $\Delta^{re}$ of \textsl{real roots} equals the $W$-orbit
of $\Pi$, which also equals $W\alpha_0$.
Denoting 
\begin{equation}
 B:=\{\gamma\in Q_+\,;\,\supp\beta\text{ is connected and }
 (\gamma,\alpha)\le 0\quad(\forall\alpha\in\Pi)\},
\end{equation}
the set of \textsl{positive imaginary roots} 
$\Delta^{im}_+$ equals $WB$.
Here
\begin{equation}
\supp\gamma:=\{\alpha\in\Pi\,;\,n_\alpha\ne0\}\text{ \  for  \ }
\gamma=\sum_{\alpha\in\Pi} n_\alpha\alpha\in\mathfrak h.
\end{equation} 

Put $\Delta^{re}_\pm:=\Delta^{re}\cap Q_{\pm}$. 
Then $\Delta^{re}=\Delta^{re}_+\sqcup 
\Delta^{re}_-$ and $WB\subset Q_+$.  
The set $\Delta$ of roots equals 
$\Delta^{re}\sqcup\Delta^{im}$ by denoting
$\Delta_-^{im}=-\Delta_+^{im}$ and 
$\Delta^{im}=\Delta_+^{im}\sqcup\Delta_-^{im}$.
The set of positive roots $\Delta_+$ equals 
$\Delta\cap Q_+=\Delta^{re}_+\sqcup WB$.
A subset $L\subset\Pi$ is called \textsl{connected}
if the decomposition $L_1\cup L_2= L$ with 
$L_1\ne\emptyset$ and $L_2\ne\emptyset$ always implies the 
existence of $v_j\in L_j$ satisfying $(v_1|v_2)\ne0$.
Note that $\alpha_0\in\supp\alpha$ for $\alpha\in\Delta^{im}$.

For $\alpha,\,\beta\in B$ with $\alpha\ne\beta$ satisfy $W\alpha\cap W\beta=\emptyset$.
Hence the set $\{\alpha_0\}\sqcup B\sqcup(-B)$ gives the complete representatives 
of $W$-orbits in $\Delta$.

There are 4 extended Dynkin diagrams $\tilde D_4$, $\tilde E_8$, $\tilde E_7$ and 
$\tilde E_6$ corresponding to connected star-shaped 
Dynkin diagrams of finite type as follows.  
Here ``finite type" means that the corresponding Weyl group is a finite group.

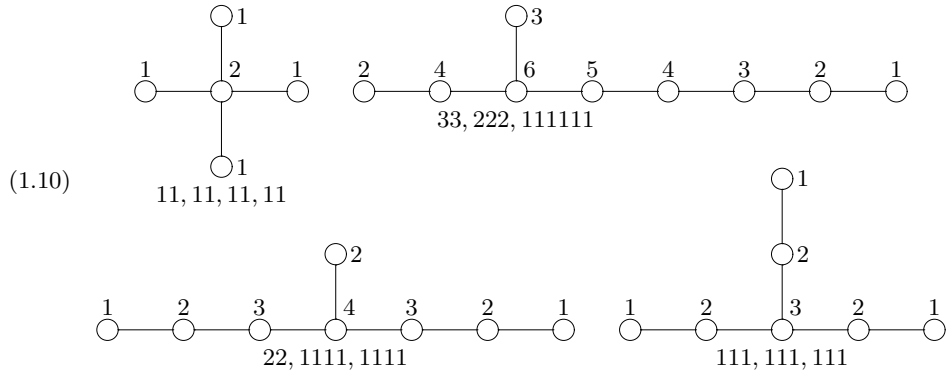
\index{Dynkin diagram}
\index{affine}
{\small
\begin{equation}\index{Dynkin diagram}\label{eq:Dynkinidx0}
\begin{gathered}
\begin{xy}
\ar@{-}               *++!D{1}  *\cir<4pt>{};
             (10,0)   *+!L+!D{2}*\cir<4pt>{}="A",
\ar@{-} "A"; (20,0)   *++!D{1}  *\cir<4pt>{},
\ar@{-} "A"; (10,-10) *++!L{1}  *\cir<4pt>{},
\ar@{-} "A"; (10,10)  *++!L{1}  *\cir<4pt>{},
\ar@{} (10,-14) *{11,11,11,11}
\end{xy}
\quad
\begin{xy}
\ar@{-}               *++!D{2}  *\cir<4pt>{};
             (10,0)   *++!D{4}  *\cir<4pt>{}="A",
\ar@{-} "A"; (20,0)   *+!L+!D{6}*\cir<4pt>{}="B",
\ar@{-} "B"; (30,0)   *++!D{5}  *\cir<4pt>{}="C",
\ar@{-} "C"; (40,0)   *++!D{4}  *\cir<4pt>{}="D",
\ar@{-} "D"; (50,0)   *++!D{3}  *\cir<4pt>{}="E",
\ar@{-} "E"; (60,0)   *++!D{2}  *\cir<4pt>{}="F",
\ar@{-} "F"; (70,0)   *++!D{1}  *\cir<4pt>{},
\ar@{-} "B"; (20,10)  *++!L{3}  *\cir<4pt>{}
\ar@{} (20,-4) *{33,222,111111}
\end{xy}
\allowdisplaybreaks\\[-.8cm]
\begin{xy}
\ar@{-}               *++!D{1}  *\cir<4pt>{};
             (10,0)   *++!D{2}  *\cir<4pt>{}="A",
\ar@{-} "A"; (20,0)   *++!D{3}  *\cir<4pt>{}="B",
\ar@{-} "B"; (30,0)   *+!L+!D{4}*\cir<4pt>{}="C",
\ar@{-} "C"; (40,0)   *++!D{3}  *\cir<4pt>{}="D",
\ar@{-} "D"; (50,0)  *++!D{2}   *\cir<4pt>{}="E",
\ar@{-} "E"; (60,0)  *++!D{1}   *\cir<4pt>{},
\ar@{-} "C"; (30,10)  *++!L{2}  *\cir<4pt>{},
\ar@{} (30,-4) *{22,1111,1111}
\end{xy}
\quad
\begin{xy}
\ar@{-}               *++!D{1}  *\cir<4pt>{};
             (10,0)   *++!D{2}  *\cir<4pt>{}="O",
\ar@{-} "O"; (20,0)   *+!L+!D{3}*\cir<4pt>{}="A",
\ar@{-} "A"; (30,0)   *++!D{2}  *\cir<4pt>{}="B",
\ar@{-} "B"; (40,0)   *++!D{1}  *\cir<4pt>{}
\ar@{-} "A"; (20,10)  *++!L{2}  *\cir<4pt>{}="C",
\ar@{-} "C"; (20,20)  *++!L{1}  *\cir<4pt>{},
\ar@{} (20,-4) *{111,111,111}
\end{xy}
\end{gathered}\end{equation}}
\end{rem}
\smallskip

We define elements $\alpha$ of $Q_+$ whose supports are given by diagrams. 
Here $\alpha$ are defined by the above diagrams with numbers attached to simple roots
which give the coefficients in the expression \eqref{eq:root}. They are also 
distinguished by the symbols $\tilde D_4$, etc.

Then the elements corresponding to the above 4 diagrams are in $B$
with $(\alpha|\alpha)=0$ and moreover it is known that any $\gamma\in B$ with 
$(\gamma|\gamma)=0$ is a positive integer multiple of one of these 4 
roots (cf.~Remark~\ref{rem:fund} (i)).

\begin{rem} 
For an integer $N$, we put
\begin{equation}
   B_N:=\{\alpha\in B\mid (\alpha|\alpha)=N\}.
\end{equation}
Then
\[B_N\ne\emptyset\ \Leftrightarrow N\in-2\mathbb Z_{\ge0}.\]
The elements of $B_0$ are multiples of the above 4 indivisible elements.
Moreover it is proved by \cite{Orims, Ow} that $B_N$ are finite sets for $N<0$
(cf.~Remark~\ref{rem:finite}).
\end{rem}

\begin{lem}[{\cite[Lemma 7.2]{Orims}}] \label{lem:root}
For the expression \eqref{eq:root} of 
$\alpha\in\Delta_+$ with $\supp\alpha\varsupsetneqq\{\alpha_0\}$, we have
\begin{equation}
\begin{gathered}
  n\ge n_{j,1}\ge n_{j,2}\ge n_{j,3}\ge \cdots\qquad
  (j=0,1,\ldots),\\
  n\le \sum n_{j,1} - \max\{n_{j,1},n_{j,2},\ldots\}.\label{eq:root0}
\end{gathered}
\end{equation}
Moreover $\alpha\in Q_+\setminus\{0\}$ is an element of $B$ if and only if
\begin{align}
  2n_{j,\nu}&\le n_{j,\nu-1}+n_{j,\nu+1}
  \quad(n_{j,0}:=n,\ j=0,1,\ldots,\ \nu=1,2,\ldots),\label{eq:neg1}\\
  2n&\le n_{0,1}+n_{1,1}+n_{2,1}+\cdots+n_{p-1,1}\label{eq:neg0}.
\end{align}
Note that these two conditions correspond to 
$(\alpha|\alpha_{j,\nu})\le0$ and $(\alpha|\alpha_0)\le0$.
\end{lem}

In this note we give an effective algorithm to calculate $B_N$
implemented in a computer program.
The root is transformed into a tuple of partitions of an integer, which
has direct applications to several fields of mathematics and mathematical physics.
We also explain an algorithm to get $\Delta_+$.

Under the expression \eqref{eq:root} of $\alpha\in\Delta_+$ with 
$\supp\alpha\supset\{\alpha_0\}$, put 
\begin{equation}\label{eq:n2m}
  m_{j,\nu}=n_{j,\nu-1}-n_{j,\nu}\quad(j=1,\dots,p-1,\ \nu=1,\dots,n_p)
.
\end{equation}
Then 
\begin{align}
   n_{j,k}&=m_{j,k+1}+\cdots+m_{j,n_j}\quad(k=1,\dots,n_j{-}1,\ j=0,\dots,p{-}1),\label{eq:m2n}\\
         n&=m_{j,1}+\cdots+m_{j,n_j}\quad(j=0,\dots,p-1)\label{eq:r2t}
\end{align}
with $m_{j,\nu}\in\mathbb Z_{\ge0}$ and we have a tuple of $p$ partitions \eqref{eq:r2t} of $n$.
Moreover we have
\begin{equation}
 (\alpha|\alpha)=2n^2-\sum_{j=0}^{p-1}\Bigl(n^2-\sum_{\nu=1}^{n_j}m_{j,\nu}^2\Bigr).
\end{equation}

The action of $W$ on $Q$ can be translated to that of 
tuples of integers :
\begin{align}
 s_{j,\nu}&: m_{j,\nu-1}\leftrightarrow m_{j,\nu}\quad(j=0,1,\ldots,\ \nu=1,2,\ldots),
\label{def:sinu}\\
 s_0&: m_{j,1}\mapsto m_{j,1}+\Bigl(2n-\sum_{j=0}^{p-1}(n -m_{j,1})\Bigr)\label{def:s0}
\quad(j=0,1,\ldots). 
\end{align}

Here we note that, for a fixed $j$, reflections $s_{j,\nu}$ $(\nu=1,2,\ldots)$ generate
permutations of indices $i$ of $m_{j,i}$ with allowing $j=p$ or  $\nu=n_j$ 
and in particular $p$ changes into $p+1$ by $s_{p,1}$.

\medskip
\begin{rem}
Under the notation in the next section,
a tuple of $p$ partitions $\mathbf m$ corresponds to an element $\alpha\in B$ if and only if
\eqref{eq:tuple}, \eqref{eq:monotone} and \eqref{eq:fund} are satisfied.  
Note that $n_{j,1}=n-m_{j,1}$ and therefore these three conditions correspond
to \eqref{eq:r2t}, \eqref{eq:neg1} and \eqref{eq:neg0}, respectively.
\end{rem}

\section{Fundamental tuples of partitions}
Let $\mathbf m$ be a tuple of $p$ partitions of a positive integer $n$, namely,
\begin{equation}\label{eq:tuple}
\begin{aligned}
   \mathbf m&=(\mathbf m_0,\mathbf m_1,\dots,\mathbf m_{p-1}),\\
   \mathbf m_j&=(m_{j,1},\dots,m_{j,n_j})&&(j=0,\dots,p-1),\\
            n&=m_{j,1}+\cdots+m_{j,n_j}&&(m_{j,\nu}\in\mathbb Z_{>0}).
\end{aligned}\end{equation}
For simplicity, we may assume $\mathbf m_j$ are {\sl monotone}, namely,
\begin{align}\label{eq:monotone}
    m_{j,1}\ge m_{j,2}\ge \cdots \ge m_{j,n_j}>0
\end{align}
and {\sl non-trivial}, namely,
\begin{align}\label{eq:nontrivial}
   n_j > 1
\end{align}
for $j=0,\dots,p-1$. We write as follows.
\begin{align}\begin{split}
   \mathbf m_j&=m_{j,1}m_{j,2}\cdots m_{p-1,n_{p-1}},\quad m_{j\nu}=m_{j,\nu},\\
   \mathbf m&=\mathbf m_0;\mathbf m_1;\cdots;\mathbf m_{p-1,n_{p-1}}\\
            &=m_{0,1}m_{0,2}\cdots m_{0,n_{0}};m_{1,1}\cdots m_{1,n_1};\cdots:
              m_{p-1,1}\cdots m_{p-1,n_{p-1}},\\
   m^k&=\overbrace{m\cdots m}^k.
\end{split}\end{align}
\begin{defn}
For a tuple $\mathbf m$ of $p$ partitions of $n$, we define
\begin{align}
   \ord\mathbf m&:=n=m_{j,1}+\cdots+ m_{j,n_j},\\
   \codim\mathbf m_j &:=n^2-\sum_{\nu=1}^{n_j}m_{j,\nu}^2,\label{def:codim}\\
   \idx\mathbf m&:=2n^2 - \sum_{j=0}^{p-1}\codim\mathbf m_j
   \label{def:idx}
\end{align}
and we call $\idx\mathbf m$ the index of rigidity of $\mathbf m$. 
The tuple $\mathbf m$ is {\sl basic} when
\begin{align}
  2n - \sum_{j=0}^{p-1}\Bigl(n-m_{j,1}\Bigr)\le 0\label{eq:fund}
\end{align}
and $\mathbf m$ is {\sl indivisible}, namely,
there exists no integer $k$ with $k>1$ which satisfies
$m_{j,\nu}\in k\mathbb Z\quad (j=0,\dots,p-1,\ \nu=1,\dots,n_j)$.
The tuple $\mathbf m$ is {\sl fundamental} if $\mathbf m$ is basic
or $\mathbf m=k\mathbf m'$ with  $k\in\mathbb Z_{>0}$ and a basic tuple $\mathbf m'$ 
satisfying $\idx\mathbf m'\ne0$. 
Here we remark $\idx k\mathbf m'=k^2\idx\mathbf m'$. 
\end{defn}
Since $\codim\mathbf m_j$ are even integers, so is $\idx\mathbf m$.
Hence the relation 
\begin{align}
  &\Bigr(\sum_{j=0}^{p-1}(\mathrm{ord}\,
  \mathbf m-m_{j,1})-2\cdot\mathrm{ord}\,\mathbf m\Bigl)\cdot\mathrm{ord}\,\mathbf m+
  \sum_{j=0}^{p-1}\Bigl(\sum_{\nu=1}^{n_j}(m_{j,1}-m_{j,\nu})m_{j,\nu}\Bigr)
  =-\mathrm{idx}\,\mathbf m\label{eq:identity}
\end{align}
means $p\ge 3$ and 
$-\idx\mathbf m$ is a nonnegative even integer for any fundamental tuple $\mathbf m$.
\begin{rem}
The number $\codim\mathbf m_j$ is the codimension of a conjugacy classes
of diagonalizable matrices in the space of matrices of size $n$ whose multiplicities of eigenvalues
are given by $m_{j,1},\ldots,m_{j,n_j}$ (cf.~\cite[\S3]{Orims} for general matrices).
\end{rem}

\medskip
To obtain fundamental tuples $\mathbf m$ with a given index of rigidity,
the following theorem is essential.
\begin{thm}\cite[Proposition 7.13]{Ow}\label{thm:main} \ 
If $\mathbf m$ is fundamental, then
\begin{align}
\mathrm{ord}\,\mathbf  m&\le 3|\mathrm{idx}\,\mathbf m|+6,\label{eq:bineq}\\
p\ge3 &\Rightarrow \mathrm{ord}\,\mathbf m\le |\mathrm{idx}\,\mathbf m|+2,\label{eq:b4ineq}\\
p&\le \tfrac12|\mathrm{idx}\,\mathbf m|+4.\label{eq:pineq}
\end{align}
\end{thm}
\begin{cor}\label{cor3}
 If $\mathbf m$ is fundamental and 
$\mathrm{ord}\,\mathbf m>|\mathrm{idx}\,\mathbf m|+2$, 
then $p=3$ and
\begin{equation}\label{eq:3pts}
m_{0,1}+m_{1,1}+m_{2,1}=\mathrm{ord}\,\mathbf m.
\end{equation}
\end{cor}
Here we may assume $\mathbf m$ is \textit{ordered}, namely,
\begin{align}\label{eq:ordered}
\begin{split}
& \exists N\in\mathbb Z_{>0} \text{ with }
m_{j,\nu}= m_{j+1,\nu}\quad(0\le \nu<N) \text{ \ and \ }m_{j,N}\ne m_{j+1,N}
\\&\quad \Rightarrow \ m_{j,\nu}>m_{j,\nu+1}
\end{split}
\end{align}
for $j=0,\dots,p-2$. 
\begin{rem}\label{rem:finite}
It follows from Theorem~\ref{thm:main} that there exist
only finite fundamental ordered tuples with a negative index of rigidity and 
the above corollary follows from \eqref{eq:identity}.
\end{rem}
\begin{exmp}\label{ex:ind0-2}
(i) \ \cite[Example 7.14]{Ow} \ 
For $m\in\mathbb Z_{>0}$, we define fundamental tuples 
\begin{equation}
\begin{aligned}\label{eq:Qsp}
 &D_4^{(m)}: m^2,m^2,m^2,m(m-1)1&& 
 \quad E_6^{(m)}: m^3,m^3,m^2(m-1)1\\ 
 &E_7^{(m)}: (2m)^2,m^4,m^3(m-1)1&& 
 \quad E_8^{(m)}: (3m)^2,(2m)^3,m^5(m-1)1
\end{aligned}
\end{equation}
with orders $2m$, $3m$, $4m$ and $6m$, 
respectively, and index of rigidity $2-2m$.
We note that $E_8^{(m)}$, $D_4^{(m)}$ and 
$11,11,11,\cdots$ attain the equalities in 
\eqref{eq:bineq}, \eqref{eq:b4ineq} and \eqref{eq:pineq}, respectively.

\smallskip
(ii) \ \cite[Corollary 6.3]{Ko} \ 
There are  4 basic tuples $\mathbf m$ with $\idx\mathbf m=0$ : 
\[D_4^{(1)}=1^2,1^2,1^2,1^2\qquad E_6^{(1)}=1^3,1^3,1^3\qquad E_7^{(1)}=2^2,1^4,1^4
 \qquad E_8^{(1)}=3^2,2^3,1^6\]

\smallskip
(iii) \ \cite[Proposition 8]{Orims} \ 
There are 13 fundamental tuples $\mathbf m$ with $\idx\mathbf m=-2$ : 
\begin{align*}
  &1^2,1^2,1^2,1^2,1^2 &&  21,21,1^3,1^3   && 2^2,2^2,2^2,21^2 && 3^1,2^2,2^2,1^4\\
  &  21^2,1^4,1^4      && 32,1^5,1^5       && 2^21,2^21,1^5    &&  3^2,2^21^2,1^6\\
  & 2^3,2^3,2^21^2     && 4^2,2^4,2^3 1^2  &&  4^2,2^4,2^31^2  &&  4^2,3^22, 1^8\\
  &6^2,4^3,2^51^2
\end{align*}
\end{exmp}
\begin{rem}\label{rem:fund}
(i) \ 
Let $\mathbf m$ be  an ordered  fundamental tuple with $\idx\mathbf m=0$. 
Then \eqref{eq:fund} and \eqref{eq:identity} show that $m_{j,\nu}=\frac{\ord\mathbf m}{n_j}$
and $\frac{\ord\mathbf m}{n_0}+\cdots+\frac{\ord\mathbf m}{n_{p-1}}=(p-1)\cdot\ord \mathbf m$. 
It is easy to conclude that $(n_0,n_2,\ldots,n_{p-1})$ equals $(2,2,2,2)$ or $(3,3,3)$ or
$(2,4,4)$ or $(2,3,6)$ and $\ord\mathbf m$ is the least common multiples of $n_0,\dots,n_{p-1}$
and we have the list in Example~\ref{ex:ind0-2} (ii).

(ii) \ 
The fundamental tuples $\mathbf m$ with $\idx\mathbf m=-4$ or $-6$ are given 
in \cite[\S13.1.2,\ \S13.1.4]{Ow}.  But some fundamental tuples are missing.
Namely 
$3^22,3^22,2^4$ should be added in the case $\idx\mathbf m=-4$.
Moreover $3^22,3^22,2^31^2$ and $3^22,3^21^2,2^4$ should be added 
in the case $\idx\mathbf m=-6$.
The fundamental tuples $\mathbf m$ with $\idx\mathbf m\le -8$ are obtained by a computer program
\cite{ORisa} using the computer algebra \texttt{Risa/Asir} whose algorithm will be explained in the next section.

(iii) \ 
In an expression of a partition, the numbers larger than 9 sometimes expressed 
by alphabets.
For example, $10,11,12,\ldots,16,17,\ldots$ are expressed by 
$\texttt{a},\texttt{b},\texttt{c},\ldots,\texttt{f},\texttt{g},\ldots$, 
respectively.  $11$ may be expressed by $(11)$ and $\texttt{bb}9=(11)(11)9$.
\end{rem}
\section{Algorithm obtaining fundamental tuples}
We want to obtain all the fundamental tuples $\mathbf m$ \eqref{eq:tuple} with a given index of rigidity.
We assume \eqref{eq:monotone} and \eqref{eq:ordered}.
Owing to \eqref{eq:bineq} and \eqref{eq:pineq}, there are finite tuples of partitions to be 
checked the conditions \eqref{def:idx} and \eqref{eq:fund} 
which assure that they are fundamental.

For example, consider the case  $\idx\mathbf m=-8$.
Then $\ord\mathbf m\le 30$ and $p\le 8$. 
When $\ord\mathbf m=30$, it follows from \cite{Ow} that $\mathbf m$ equals $E_8^{(8)}$.
Suppose $\ord\mathbf m=29$. Then \eqref{eq:pineq} assures $p=3$.
Since the number of divisions $p(29)=4565$, we have to check about $4565^3/6>10^{10}$ cases
of ordered triplets of partitions of 29.  This is quite heavy.
As a result, there is no tuple satisfying these condition. 

In this section we will show and explain a computer program which gives the fundamental tuples 
for a fixed index of rigidity.  It gives 116 ordered fundamental tuples $\mathbf m$
with $\idx\mathbf m=-8$ within 0.05 second under Windows 11 with CPU i7-10700.
In the same way, it takes about 1 second (resp.\ 3 minutes) to get 884 (resp.~40617) fundamental tuples $\mathbf m$
with $\idx\mathbf m=-18$ (resp.~$-50$).  Note that $p(3\cdot 50+5)=66493182097 \fallingdotseq 6.65\times 10^{10}$.

\medskip
We use a free computer algebra \texttt{Risa/Asir} \cite{Risa} and
a partition is expressed by a list of integers and a tuple is a list of lists of integers, namely,
the tuple $\mathbf m$ is represented by a list of lists \text{V} : 
\begin{align}
\begin{split}
  \mathbf m&=\bigl[[m_{0,1},\dots,m_{0,n_0}],[m_{1,1},\dots,m_{1,n_1}],\dots,
   [m_{p{-}1,1},\dots,m_{p{-}1,n_{p{-}1}}]\bigr],\\
  \texttt{V[}j\texttt{]}&\texttt{[}\nu-1\texttt{]}
    =V_{j,\nu-1}=m_{j,\nu}\quad(j=0,1,\dots,p-1, \ \nu=1,2,\ldots,n_j)
\end{split}
\end{align}
under the notation in the previous section. 
We assume $\mathbf m$ is monotone and ordered.
In our program we also use symbols
\begin{align}
\begin{split}
    \texttt{D}&=\mathrm{ord}\,\mathbf m=n,\\
    \texttt{Idx}&=|\mathrm{idx}|,\\
    \texttt{S}&:=\sum_{j=0}^p(\mathrm{ord}\,\mathbf m-m_{j,1})-2\cdot\mathrm{ord}\,\mathbf m,\\
    \texttt{SS}&:=\sum_{j=0}^p\Bigl(\sum_{\nu=1}^{n_j}(m_{j,1}-m_{j,\nu})m_{j,\nu}\Bigr).
\end{split}
\end{align}
Then the the tuple $\mathbf m$ is fundamental if and only if 
\begin{equation}\label{eq:check}
  \texttt{S}\ge0 \text{ \ and \ } \texttt{S}+\texttt{SS}=\Idx
\end{equation}
together with the condition that $\mathbf m$ is indivisible when $\Idx=0$.

\bigskip
Now we explain the function \texttt{spbasic(\,)} in a library \texttt{os\_muldif.rr}  
\cite{ORisa}.

\smallskip
\noindent
\texttt{spbasic(\Idx,\D}\verb/|/\texttt{str=1,pt=[$k$,$\ell$])}\\
:: Returns the list of ordered fundamental tuples with rank and index\ \Idx.
\smallskip

If \Idx\ equals 0, indivisible fundamental tuples are return.

If \D\ equals 0, the list of all the fundamental tuples with index \Idx\ are returned.

Here the fundamental tuples are given in the form of ordered tuples.

If the option \texttt{pt=[$k$,$\ell$]} is indicated, then $k\le \# \text{partitions (in a tuple)} 
\le\ell$

The option \texttt{pt=[$k$,$k$]} may be indicated by \texttt{pt=$k$}.

If the option \texttt{str=1} is indicated, tuples of partition are given by 
strings of numbers as in Example~\ref{ex:ind0-2}.  These two expressions of tuples of 
partition can be transformed to each other by \texttt{s2sp(\,)} in \texttt{os\_muldif.rr}. 

\bigskip
A partition is expressed by a sequence of non-increasing positive integers.
Put larger partitions first in the lexicographic order of sequences of numbers and 
put larger tuples of partitions first in the lexicographic order.
For example,  
\[\texttt{aa431}>\texttt{a99}>\texttt{a981}\qquad
 \texttt{541,3331} > \texttt{532,433}\]

The function \texttt{spbasic(\Idx,\D)} checks whether each tuple of partitions satisfies
\eqref{eq:check} according to the above order and returns all the required tuples.
Note that Theorem~\ref{thm:main} assures that  
there is no required tuple if $\D>3\cdot|\Idx|+6$ or $\Idx>0$ or $\Idx$ is odd and 
moreover that the number $p$ of partitions of the tuple
equals 3 if $\D>|\Idx|+2$ and moreover $3\le p\le \tfrac12|\Idx|+4$ in general.  
Hence the number of tuples to be checked is finite but it is quite large 
when $|\Idx|$ and $\D$ are large.
For a practical and effective program to perform this, 
it is quite important to avoid unnecessary checks 
as much as possible.
Note that most fundamental tuples $\mathbf m$ in this region satisfy $|\idx\mathbf m|>|\Idx|$.

\medskip
The code of \texttt{spbasic(\,)} is divided into two parts:

\smallskip
\noindent
{\bf 1. Introduction}

{\bf 1.1, 1.2, 1.3.}
Check the arguments and the options of the function and return the answer 
when it is clear, for example, in the case $\D\ge 3\cdot|\Idx|+6$.

{\bf 1.4.}
Determine an upper limit $\texttt{L}$ such that $p\le \texttt{L}$ according to
Theorem~\ref{thm:main}.

{\bf 1.5, 1.6.}
Determine an upper limit \texttt{T} such that $m_{0,1}\le \texttt{T}$
by using the condition $\texttt{SS}\ge 0$. Then $m_{j,\nu}\le\texttt{T}$. 
Set the first tuple to be checked with $m_{0,1}=\texttt{T}$.

\medskip
\noindent
{\bf 2. Loop of Main routine}

{\bf 2.1.}
If $\D>|\Idx|+2$, check \eqref{eq:3pts} and skip tuples 
without \eqref{eq:3pts}.

{\bf 2.2.}
Examine the minimal integer $\texttt{IL}$ such that $\texttt{S}\ge 0$ with $p=\texttt{IL}<\L$.

If $\IL$ does not exist,
skip some unsuitable tuples and go to the top of the loop 

{\bf 2.3.}
If there exists $\texttt{K}\le \IL$ such that 
$\texttt{SS}>|\Idx|$ with $p=\texttt{K}$, 
skip some unsuitable tuples and go to the top of Loop.

{\bf 2.4.}
Examine the minimal $p=\texttt{J}$ such that $\texttt{Ix}:=\Idx-\texttt{S}-\texttt{SS}\le 0$
with $\texttt{IL}\le \texttt{J}<\texttt{L}$.

{\bf 2.5.}
If \texttt{J} is found and $\texttt{Ix}=0$, a required tuple is obtained
and the tuple is saved in \texttt{R}.

{\bf 2.6.}
Skip some unsuitable tuples according to \texttt{Ix} and go to the top of Loop.
\bigskip

We show the code of \texttt{spbasic(\,)} with explanations inserted.

\bigskip
{\color{black}\begin{verbatim}  
def spbasic(Idx,D)
{
\end{verbatim}}
\noindent
// {\bf 1.} 
{\bf Initial setting}

\noindent
//\qquad \Idx\ $\mapsto$  \texttt{-\Idx}

\noindent
//\qquad Set options \texttt{str} (output by strings). 

\noindent
//\qquad Handle the option \texttt{pt=}.
{\color{black}\begin{verbatim}  
  Idx=-Idx;
  if((Str=getopt(str))!=1) Str=0;
  Tu0=Tu1=0;
  if(type(Tu=getopt(pt))==4&&length(Tu)==2&&isint(car(Tu))){
    Tu0=Tu[0];Tu1=Tu[1];
  }else if(isint(Tu)) Tu0=Tu1=Tu;
  if(Tu1<3) Tu1=0;\end{verbatim}} 
\noindent
// {\bf 1.1.}
If $\D> 3\cdot \Idx+6$, there is no solution. 
Hence we assume $\D \le 3\cdot \Idx+6$.
{\color{black}%
\begin{verbatim}
  if(!isint(Idx)||!isint(Idx/2)||Idx<0||!isint(D)||D<0||D==1
    ||D>3*Idx+6) return [];
\end{verbatim}}  
\noindent
// {\bf 1.2.}
$\D=0\ \Rightarrow\ $Handle arbitrary order.
{\color{black}\begin{verbatim}  
  if(D==0){
    for(R=[],D=3*Idx+6;D>=2;D--)
      R=append(spbasic(-Idx,D|str=Str,pt=[Tu0,Tu1]),R);
    return R;
  }
\end{verbatim}}
\noindent
// {\bf 1.3.}
Return the answer when $\Idx=0$.
{\color{black}\begin{verbatim}  
  if(!Idx){
    R=0;
    if(D==2&&Tu0<5&&Tu1>3) R="11,11,11,11";
    if(Tu0<4&&Tu1>2){
      if(D==3) R="111,111,111";
      if(D==4) R="22,1111,1111";
      if(D==6) R="33,222,111111";
    }
    if(!R) return [];
    return [(Str==1)?R:s2sp(R)];
  }
\end{verbatim}}  
\noindent
// {\bf 1.4.}
Set \texttt{L} an upper bound  of the number of partitions of a tuple. 

\noindent
//\quad $\D>\Idx+2\qquad\ \Rightarrow\ $ \#  partitions of a tuple $=3$ ($\texttt{L}=3$).

\noindent
//\qquad$\D=3\cdot\Idx+6\ \!\Rightarrow\ $ Returns the answer.

\noindent
//\quad $\D=2\qquad\qquad\quad\ \Rightarrow\ p=\Idx/2+4$

\noindent
//\quad Check the option \texttt{pt=}.
{\color{black}\begin{verbatim}  
  if(D>Idx+2){
    L=3;
    if(D==3*Idx+6){
      R=[[D/2,D/2],[D/3,D/3,D/3],[D/6,D/6,D/6,D/6,D/6,D/6-1,1]];
      return [(Str==1)?s2sp(R):R];
    }
  }else L=Idx/2+4;
  if(L>Tu1) L=Tu1;
\end{verbatim}}
\noindent
// {\bf 1.5.}
Find a positive number \texttt{T} satisfying $V_{00}\le \texttt{T}$ \ 
($\Rightarrow\, V_{j,\nu}\le\texttt{T}$).

\noindent
//\quad Note that $V_{00}<\D$ and $K\cdot(V_{00}-\texttt{K})\le\Idx$ with 
$\texttt{K}=\D\%V_{00}$. 

\noindent
//\quad Here $\D\%V_{00}$ is a minimal non-negative integer in $\D+V_{00}\mathbb Z$.

{\color{black}\begin{verbatim}  
  for(T=D-1;T>1;T--){
    K=D%T;
    if((T-K)*K<=Idx) break;
  }
\end{verbatim}}
\noindent
// {\bf 1.6.}
Define the initial \texttt{V} by the top tuple with $V_{00}=\texttt{T}$.

\noindent
//\quad \texttt{NP[$m$]} : the top partition of the integer $\D$ by positive integers $\le m$. 

\noindent
//\quad \texttt{NP[$m$]}$=$\texttt{[$m,\dots,m,m'$]} with \texttt{length(NP)}$=\D$ and 
$m\ge m'>0$.

\noindent
//\quad \texttt{FS}: check \eqref{eq:3pts} ?
{\color{black}\begin{verbatim}  
  V=newvect(L);
  NP=newvect(T+1);
  for(I=1;I<=T;I++) NP[I]=nextpart(D|max=I);
  for(I=0;I<L;I++) V[I]=NP[T];
  S1=NP[1];
  FS=(D>Idx+2 ||(L==3&&D>Idx))?1:0;
\end{verbatim}}

\medskip
\noindent
// {\bf 2.} {\bf Loop of Main routine }

\noindent
//\quad \Idx\ (index) ,\ \D\ (order), \L\ (maximal number of partitions of a tuple) are given.

\noindent
//\quad \texttt{V} (initial tuple) is also given.  The required tuples $\le $ \texttt{V}. 
{\color{black}\begin{verbatim}  
  for(R=[];;){
\end{verbatim}}
\noindent
// {\bf 2.1.}
Case  $\D>\Idx+2$ etc. : examine $V_{00}$, $V_{10}$, $V_{20}$ \quad ($V_{00}\ge V_{10}\ge V_{20}$)

\noindent
//\quad Check the condition  $V_{00}+V_{10}+V_{20}=\D$.

\noindent
//\quad $3V_{00}<\D\ \Rightarrow\ $ no more required tuples.
{\color{black}\begin{verbatim}  
    if(FS){
      if(3*car(V[0])<D) break;
\end{verbatim}}  
\noindent
// {\bf 2.1.1.}
$V_{00}+V_{10}\ge\D\ \ \Rightarrow \ \ $ change $V_{10}\mapsto \D - V_{00}-1$ if possible.
{\color{black}\begin{verbatim}  
      if(car(V[0])+car(V[1]) >= D && (T=D-car(V[0])-1) > 0)
        V[1]=V[2]=NP[T];
\end{verbatim}}  
\noindent
// {\bf 2.1.2.}
Want to have the equality $\texttt{S}\,:=\D-V_{00}-V_{10}-V_{20} = 0$.
{\color{black}\begin{verbatim}  
      S=D-car(V[0])-car(V[1])-car(V[2]);
\end{verbatim}}
\noindent
// {\bf 2.1.3.}
$V_{00}+2\cdot V_{10}<\D$ or $(\texttt{S}\ \ne 0,\ V_{10}=1)\ \Rightarrow\ $change $ V_{0},\ldots$
$\Rightarrow$ Top of Loop.
{\color{black}\begin{verbatim}  
      if(car(V[0])+2*car(V[1])<D || (car(V[1])==1 && S) ){
        if(car(V[0])==1) break;
        V[0]=V[1]=V[2]=nextpart(V[0]);
        continue;
      }else if(S<0) V[2]=NP[car(V[2])+S];
\end{verbatim}}  
\noindent
// {\bf 2.1.4.}
$S>0$ or $V_{20} \le S_2\ \ \Rightarrow\ $ change \texttt{V[1]} with 
\texttt{V[2]}\ 
$\Rightarrow$\ Top of Loop.
{\color{black}\begin{verbatim}  
      if(S>0||var(V[2])+S<1){
        V[1]=V[2]=nextpart(V[1]);
        continue;
\end{verbatim}}  
\noindent
// {\bf 2.1.5.}
$V_{20}$ can be changed so that $V_{00}+V_{10}+V_{20}=\D$.
{\color{black}\begin{verbatim}
      }else if(S<0) V[2]=NP[car(V[2])+S];
    }
\end{verbatim}}

\medskip
\noindent
// {\bf 2.2.}
Find minimal \IL\ satisfying $(\D{-}V_{00}){+\cdots +}(\D{-}V_{\IL0}){-}2D\ge 0$ with $\IL<\L$.
{\color{black}\begin{verbatim}  
    for(S=-2*D,IL=0;IL<L;IL++){
      S+=D-car(V[IL]);
      if(S>=0) break;
    }
\end{verbatim}}

\noindent
// {\bf 2.2.1.}
$\texttt{S}=(\D-V_{00})+\cdots+(\D-V_{L{-}1,0}) - 2\D<0$. 

\noindent//\quad
$\Rightarrow$ Get minimal \texttt{K}  so that $V_{\texttt{K}0}$ 
can be reduced to attain \eqref{eq:fund}
$\Rightarrow$\ Top of Loop.
{\color{black}\begin{verbatim}  
    if(S<0){  /* reducible i.e. IL=L && S<0 */
      for(LL=L-1;LL>=0;LL--){
        if((K=car(V[LL]))+S>0){
          V[LL]=NP[K+S];
          break;
        }else{
          S+=K-1;
          V[LL]=S1;
        }
      }
      if(LL<0) break;
      for(I=LL;I<L;I++) V[I]=V[LL];
      continue;
    }
\end{verbatim}}  
\noindent
// {\bf 2.3.}
We have $\texttt{S}:=(\D-V_{00})+\cdots+(\D-V_{k,0}) - 2\D\ge0$ with $\IL\le k<\L$.

\noindent//\quad 
Check the second term \texttt{SS} of LHS of \eqref{eq:identity} which should be $\le \Idx$.
{\color{black}\begin{verbatim}  
    for(SS=K=0;K<=IL;K++){
      ST=car(V[K]);SS0=SS;
      for(I=length(V[K])-1;I>0;I--) SS+=(ST-V[K][I])*V[K][I];
      if(SS>Idx) break;
    }
\end{verbatim}}
\noindent
// {\bf 2.3.1.}
$\texttt{SS}> \Idx$ for $\texttt{V[0]},\ldots,\texttt{V[K]}$ with $\texttt{K}\le\texttt{IL}\\
//\qquad\quad \Rightarrow\ $change \texttt{V[K]}$\ \Rightarrow\ $ Top of Loop.
{\color{black}\begin{verbatim}  
    if(SS>Idx && car(V[K])!=1){
      if((W=nextcod(V[K],SS-Idx))==[]){
        for(T=car(V[K])-1;T>0;T--){
          J=D%T;
          if(SS0+J*(T-J)<=Idx) break;
        W=NP[T];
      }
      for(J=K;J<L;J++) V[J]=W;
      continue;
    }
\end{verbatim}}
\noindent
// {\bf 2.4.}
Check the tuples $\texttt{V[0]},\ldots,\texttt{V[J]}$ whether they attain $\Idx$.

\noindent//\quad
\Ix: \Idx $\,-\,$(LHS of \eqref{eq:identity} with $p=\texttt{J}$), which strictly decreases 
when \texttt{J} increases.

\noindent
//\quad \IxF: the value of \Ix\ in the case \texttt{J}${}-1$.

\noindent
//\quad Get smallest \texttt{J} with $\Ix\le 0$.
{\color{black}\begin{verbatim}  
    for(Ix=2*D^2+Idx,J=0;J<L;J++){
      IxF=Ix;
      for(Ix-=D^2,TV=V[J];TV!=[];TV=cdr(TV)) Ix+=car(TV)^2;
      if(Ix<=0) break;
    }
\end{verbatim}}
\noindent
// {\bf 2.5.}
$\Ix=0\ \Rightarrow\ $ we have a required tuple $\Rightarrow$ save it in \texttt{R}.
{\color{black}\begin{verbatim}  
    if(!Ix&&J>=IL&&J>Tu0-2){
      for(TR=[],K=J;K>=0;K--) TR=cons(V[K],TR);
      R=cons((Str==1)?s2sp(TR):TR,R);
    }
\end{verbatim}}
\noindent
// {\bf 2.6.}
Prepare the next search.

\noindent
// {\bf 2.6.1.}
$\Ix<0\ \Rightarrow\ $ check the case that \texttt{V[J-1]} should be changed.
{\color{black}\begin{verbatim}  
    if(Ix<0 && IxF-mincod(D,car(V[J]))<0) J--;
\end{verbatim}}  
\noindent
// {\bf 2.6.2.}
Cannot find \texttt{J}\ $\Rightarrow$\ 
change \texttt{V[\L-1]} in 2.6.3.
{\color{black}\begin{verbatim}  
    else if(J>=L) J=L-1;
\end{verbatim}}
\noindent
// {\bf 2.6.3.}
Change \texttt{V[J]} to the next one.
{\color{black}\begin{verbatim}
    for(I=J;I>=0&&car(V[I])==1;I--);
    if(I<0) break;
    V[I]=nextpart(V[I]);
    for(J=I+1;J<L;J++) V[J]=V[I];
  }
  return R; 
}
\end{verbatim}}

\bigskip
\rm
{\bf 3.}
Functions  used in \texttt{spbasic()}

\noindent
{\bf 3.1.} 
Get the next partition.\\
\noindent
\texttt{nextpart(V}\verb/|/{\texttt{max=N})}\\
:: get the partition that comes after \texttt{V} in descending lexicographical order.

If \texttt{V} is the last partition, \texttt{0} is returned.

If \texttt{V} is a positive integer, the top partition of the integer \texttt{V} is given.
In this case, the option \ \texttt{max=N} \ means the top partition of \texttt{V}
by integers $\le \texttt{N}$, which equals \texttt{[N,$\ldots$,N,M]} with  
$0<\texttt{M}\le\texttt{N}$.

Note that any partition $\mathbf m$ of 
\texttt{V} in the position after \texttt{nextpart(V}\verb/|/{\texttt{max=$N$)}
satisfies \ $\codim\mathbf m >\codim \texttt{nextpart(V}$\verb/|/$\texttt{max=N})$.

\medskip
We show examples executing \texttt{nextpart(\ )}.  
For the execution, we should add the prefix ``\verb|os_md.|" to 
the name of the function in the library \verb|os_muldif.rr| \cite{ORisa}.
\begin{verbatim}
[0] S=os_md.nextpart(5);
    [5]
[1] while((S = os_md.nextpart(S)) != 0) print(S);
    [4,1]
    [3,2]
    [3,1,1]
    [2,2,1]
    [2,1,1,1]
    [1,1,1,1,1]
    0
[2] os_md.nextpart(5|max=3);
    [3,2]
\end{verbatim}
{\color{black}\begin{verbatim}
def nextpart(V)
{
  if(isint(V)){
    if(V<1) return [0];
    I=V;
    if(!isint(K=getopt(max))|| K<1) K=I;
    K++;V=[];
  }else{
    if(car(V) <= 1)
      return 0;
    for(I = 0, V = reverse(V); car(V) == 1; V=cdr(V))
      I++;
    I += (K = car(V));
    V=cdr(V);
  }
  R = irem(I,--K);
  R = (R==0)?[]:[R];
  for(J = idiv(I,K); J > 0; J--)
    R = cons(K,R);
  while(V!=[]){
    R = cons(car(V), R);
    V = cdr(V);
  }
  return R;
}
\end{verbatim}}
\medskip

{\bf 3.2.} 
Get the next partition with a required property\\
\texttt{nextcod(V,D)}\\
:: get the partition after \texttt{V} whose square sum of elements $\ge$ (that of \texttt{V}${}){}+\D$.

If the required partition does not exist, \texttt{[]} is returned.
{\color{black}\begin{verbatim}
def nextcod(V,D)
{
  if(D<=0){
    if(car(V)==1) return [];
    for(D=0,TV=V;TV!=[];TV=cdr(TV)) D+=car(TV)^2;
    V=nextpart(V);
    for(TV=V;TV!=[];TV=cdr(TV)) D-=car(TV)^2;
  }
  K=length(V)-1;
  for(S=T=0;K>=0; K--){
    T+=(W=V[K]);S+=W^2;
    if(W>2&&T^2-mincod(T,W-1)>=S+D){
      R=nextpart(T|max=W-1);
      for(--K;K>=0;K--) R=cons(V[K],R);
      return R;
    }
  }
  return [];
}
\end{verbatim}}

\medskip
{\bf 3.3.} 
Get the minimal codimension of certain conjugacy classes of matrices\\
\texttt{minicod(D,T)}\\
:: returns $\codim \texttt{nextpart(D,T)}$  (cf.~\eqref{def:codim})

{\color{black}\begin{verbatim}
def mincod(D,T)
{
  K=D%T;
  return D^2-T*(D-K)-K^2;
}
\end{verbatim}}

\smallskip
We examine the effects of codes in \texttt{spbasic(\,)} eliminating unnecessary checks
and the runtime of \texttt{spbasic(\Idx,0)} is given under Windows 11 with CPU i7-10700.

\medskip
\noindent\ \ 
\begin{tabular}{|c|r|r|r|r|r|}\multicolumn{6}{c}
{\bf Effect of codes eliminating unnecessary checks}
\\ \hline
Code                   & 1.5    & 2.1     & 2.31   & 2.31${}^*$ &2.61   \\\hline
check                  & $\texttt{SS}>|\Idx|$ & $\texttt{S}\ne0$ & $\texttt{SS}>|\Idx|$
 & $\texttt{SS}>|\Idx|$ &$\texttt{Ix}<0$  \\\hline
Index of rigidity      & $-16$ & $-26$ & $-8$ &$-16$ &$-26$  \\\hline
\# fundamental tuples&647 &2889 &291 &647&2889 \\\hline
original runtime (sec)& 0.58  & 4.92   & 0.05  &0.58& 4.92  \\\hline
runtime without Code & 26.31 & 24.47   & 91.58  &32.44& 24.67  \\\hline
\end{tabular}

\smallskip
Here we remark $\texttt{Ix}<0\ \Leftrightarrow \ \texttt{S}+\texttt{SS}>|\Idx|$.

The term ``runtime without Code" in the above column 2.31$^*$ means the runtime with using 
\texttt{nextpart(\,)} in place of  \texttt{nextcod(\,)} in \textbf{2.31}.
\begin{rem}
Fundamental tuples whose index of rigidity $\ge-10$ are listed in \cite{Ofund}
and rigid tuples of order at most 8 are listed in \cite[13.2.3]{Ow}.
\end{rem}

\section{Related results}
For a given $\alpha\in\Delta$, it is easy to get the orbit $W\alpha\subset\Delta$.
Namely, the function \texttt{spgen(\,)} in \texttt{os\_muldif.rr} \cite{ORisa} gives the finite set
$\{\gamma\in W\alpha\mid \ord\gamma\le \D\}$ for $\D\in\mathbb Z_{>0}$.
Since the set is stable under the transformation by $\sigma_{j,j}\in S_\infty$, 
it is better to give the representatives of the elements of the set as monotone tuples
of partitions.

For a tuple of partitions $\mathbf m$, define $s\mathbf m$ the unique monotone tuple obtained by 
successive applications of $s_{j,\nu}$ (cf.~\eqref{def:sinu}). 

For an integer $r\ge p-1$ and 
$(j_0,\ldots,j_r)\in \mathbb Z_{>0}^{r+1}$, put
\begin{align}
\bar T_{\nu_0,\nu_1,\ldots,\nu_r}&: m_{j,\nu_j}\mapsto  m_{j,\nu_j}-\Bigl(2n-\sum_{j=0}^r(n -m_{j,\nu_j})\Bigr)\quad(j=0,\dots,r),\label{eq:red1}\\
T_{\nu_0,\nu_1,\ldots,\nu_r}&=s\circ\bar T_{\nu_0,\nu_1,\ldots,\nu_r}.
\end{align}
We note that $\bar T_{\nu_0,\nu_1,\ldots,\nu_r}$ is an involutive transformation 
corresponding to an element of $W$ and moreover 
\begin{equation}\label{eq:diford}
  \ord T_{\nu_0,\nu_1,\ldots,\nu_r}\mathbf m-\ord\mathbf m
  =\sum_{j=0}^r(\ord\mathbf m -m_{j,\nu_j})-2\ord\mathbf m.
\end{equation}
Owing to \eqref{eq:root0}, the right hand side of \eqref{eq:diford} is positive if $r\ge p$.

It is known that for a fundamental tuple $\mathbf m$, 
any monotone tuple of partitions $\mathbf m'\in W\mathbf m$ has an expression
\begin{align}\label{eq:generate}
\begin{split}
  &\mathbf m^{(0)}=\mathbf m,\ \mathbf m^{(N)}=\mathbf m',\ \ 
  \mathbf m^{(k)}=T_{j_0^{(k)},\dots,j_{r_k}^{(k)}}\mathbf m^{(k-1)},\
  \ord \mathbf m^{(k)}>\ord \mathbf m^{(k-1)}
\end{split}
\end{align}
with suitable $T_{j_0^{(k)},\dots,j_{r_k}^{(k)}}$ for $k=1,\dots,N$.

Let $\gamma\in\Delta^{re}_+$ with $\supp\gamma\ni\alpha_0$.  
Then we may assume the corresponding tuple $\mathbf m'$ of partitions is monotone and 
$\mathbf m'$ is similarly obtained by \eqref{eq:generate} with putting
$\mathbf m$ the trivial tuple of partition $1$ (or $1,1,1$).

By the sequence \eqref{eq:generate}, 
we obtain an algorithm to get tuples of partitions
corresponding to elements of $\Delta_+$, which is realized as follows. 

\medskip
\noindent
\texttt{spgen}%
$\verb|(|n\verb/|eq=1,str=1,std=/f\verb/,pt=[/k,\ell\verb/],sp=/m\verb/,basic=1)/$\\
:: Get tuples $\mathbf m$ of partitions corresponding to the elements $\Delta_+$ with 
$\supp\alpha\ni\alpha_0$.

\begin{itemize}
\item
If $n\in\mathbb Z_{>0}$, rigid tuples $\mathbf m$  with $\ord\mathbf m\le n$
are obtained and  
\begin{itemize}
\item
\verb|eq=1|\phantom{AAa} : $\ord\mathbf m=\D$
\item
\verb|str=1|\phantom{AA} : returns in string form
\item
\verb|pt=[|$k,\ell$\verb|]| : the number of the partition is in $[k,\ell]$. 
\texttt{pt=$k$} means $k=\ell$
\item
\texttt{std=$f$}\phantom{Aa} : the multiplicities are arranged from larger ones at each partition
and the partitions are arranged from smaller (resp.~larger) tuples if $f=1$
(resp.~$f=-1$).
\end{itemize}
\item
If  $n\in\mathbb Z_{\le0}$, the fundamental tuples $\mathbf m$ with $\idx\mathbf m=n$ are obtained.
Options \verb|str=1| and \verb|pt=| are valid.
\item
The option \verb|sp=|$\mathbf m$  returns tuples in $W\mathbf m$ obtained by the 
direction with higher rank. 
If \verb|basic=1| is indicated, the direction is not specified.
Options \verb|st=1|, \verb|std=| and \verb|pt=| are valid.
\end{itemize}
\begin{verbatim}
[0] os_md.spgen(4|eq=1,pt=4);
 [[[2,2],[2,2],[2,2],[3,1]],[[2,2],[3,1],[3,1],[2,1,1]]]
[1] ltov(os_md.spgen(4|eq=1,pt=4,str=1));
 [ 22,22,22,31 22,31,31,211 ]
[2] ltov(os_md.spgen(4|eq=1,pt=4,str=1,std=-1));
 [ 31,22,22,22 31,31,22,211 ]
[3] ltov(os_md.msort(os_md.spgen(4|eq=1,pt=4,str=1,std=-1),[-1,0]));
 [ 31,31,22,211 31,22,22,22 ]
[4] ltov(os_md.spgen(-2|pt=4,str=1));
 [ 21,21,111,111 22,22,22,211 31,22,22,1111 ]
[5] ltov(os_md.spgen(4|eq=1,sp="11,11,11,11",str=1,std=-1));
 [ 31,31,22,1111 31,22,22,211 31,31,211,211 31,31,31,31,22 ]
\end{verbatim}
Here the header \texttt{os\_md.} is required 
for a function in \texttt{os\_muldif.rr}.

\texttt{[0]} list of rigid tuples $\mathbf m$ with $p=4$ and $\ord\mathbf m=4$

\texttt{[1]} the same list as above in the string form

\texttt{[2]} the same list as above in the string form with the inverse lexicographic order

\texttt{[3]} the same list as above arranged by tuples

\texttt{[4]} list of fundamental tuples $\mathbf m$ with $\idx\mathbf m=-2$ and $p=4$

\texttt{[5]} list of tuples $\mathbf m$ with $\ord\mathbf m=4$ in the orbit containing 11,11,11,11

\bigskip
A relation between a star-shaped Kac-Moody root system and Deligne-Simpson problem of Fuchsian
ordinary differential equations was first clarified by \cite{CB}.
For a Fuchsian differential equation 
on the Riemann sphere with $p$ singular points
\begin{align*}
 \mathcal M: 
\frac{du}{dx}=\Bigl(\sum_{j=1}^{p-1}\frac{A_j}{x-a_j}\Bigr)u\ \ \text{with } A_1+\cdots+A_{p-1}+A_p=0
\quad(A_j\in M(n,\mathbb C))
,
\end{align*}
the spectral type of $\mathcal M$ is a tuple of $p$ partitions of $n$ defined by 
the multiplicities $m_{j,\nu}$ of characteristic exponents $\lambda_{j,\nu}$,
namely the eigenvalues of residue matrices
$A_j$ at the singular points $a_j$ ($a_p=\infty$) when $A_j$ are diagonalizable.  
Here $n$ is the rank of $\mathcal M$.

Deligne-Simpson problem is to find a necessary and sufficient condition which assures 
the existence of an irreducible tuple $(A_1,\dots,A_p)$ of matrices with $A_1+\cdots+A_p=0$ 
in terms of conjugacy classes of $A_\nu$  ($\nu=1,\dots,p$). 
The condition $\trace A_1+\cdots+\,\trace A_p=0$ is clearly necessary.
The tuple said to be  
irreducible if $A_j V\subset V$ $(j=1,\dots,p)$ for a subspace $V$ of $\mathbb C^n$
means $V=\{0\}$ or $V=\mathbb C^n$.
Note that if $A_j$ is diagonalizable, the conjugacy class is determined by its eigenvalues
$\lambda_{j,\nu}$ and their multiplicities $m_{j,\nu}$  ($\nu=1,\dots,n_j$).

There exists an irreducible Fuchsian differential equation with a given spectral type 
if and only if the tuple $\mathbf m$ corresponds to a positive root $\alpha$ of 
the star-shaped Kac-Moody root system with the condition that $\alpha_0\in \supp\alpha$ 
and moreover $(\alpha|\alpha)\ne 0$ or $\alpha$ is indivisible
(cf.~\cite{CB, Ow}).  
Under this condition, the irreducible equation exists if 
$\lambda_{j,\nu}$ are generic under the condition $\sum_{j,\nu}m_{j,\nu}\lambda_{j,\nu}=0$.

Then $2-(\alpha|\alpha)$ gives the number of accessory 
parameters of the equation.  
Here the accessory parameters correspond the 
moduli describing the simultaneous conjugacy classes of 
$(A_1,\dots,A_j)$ for given $\lambda_{j,\nu}$ and $m_{j,\nu}$. 
The equation is called {\sl rigid} if it has no accessory parameter and
it means the spectral type corresponds to a positive real root. 

To get this result, transformations of the equations called middle convolutions and additions 
are essential, which induce transformations of spectral types 
corresponding to the action of elements of $W$ on the Kac-Moody root 
system

As in the case of Fuchsian system $\mathcal M$, 
\cite{Ow} gets a similar result for single linear ordinary differential equations 
with regular singularities in $\mathbb P^1$.
In fact, middle convolution, addition, spectral type are similarly defined and
moreover various problems including 
construction of the equations, integral representation of their solutions, their 
connection problem etc. are studied.
The star-shaped Kac-Moody root system and the action of the Weyl group is essential
in the study.

\medskip
Let $\alpha\in\Delta_+$ with $\supp\alpha\not\ni \alpha_0$. 
Then $\alpha$ correspond to a classical 
root of type $A$, namely, there exist $j_0\in Z_{\ge0}$ and non-negative integers 
$\nu_0<\nu_1$ such that
\begin{equation}
  m_{j,\nu}=\begin{cases}
            0& (j,\nu)\ne (j_0,\nu_0),\ (j_0,\nu_1),\\
            (-1)^{\epsilon} &(j=j_0,\ \nu=\nu_\epsilon,\ \epsilon=0,\,1).
			\end{cases}
\end{equation}
This tuple $\mathbf m$ has a special meaning for reducibility of 
rigid Fuchsian differential equations 
(cf.~\cite[Type 3 in Definition~2.4]{Ored}).

\bigskip
We explain a function \texttt{chksp(\,)} in \texttt{os\_muldif.rr} which analyzes a tuple 
of partitions of a positive integer. Here the tuple may not be monotone nor ordered.
It also analyzes generalized Riemann schemes but here we do not explain it
(cf.~\cite[\texttt{os\_muldifeg.pdf}]{ORisa}).

\medskip
\noindent
\texttt{chkspt}%
$\verb|(|\mathbf m\verb/|opt=/t\verb|)|$\\ 
:: Check a tuple  $\mathbf m$ of partitions (spectral type) or 

\qquad\,a generalized Riemann scheme 
$\{[\lambda_{j,\nu}]_{m_{j,\nu}}\}$\\
\quad Return : 
\texttt{[$pts,\, ord,\, idx,\, fuchs,\, rod,\, redsp, fspt$]}

\qquad\ \ \  
$-1$ means the tuple is not correct

\qquad\ \ \ 
$0$ means that the tuple is not in $\Delta_+$
\begin{itemize}
\item[$pts\phantom{sp}$] : number $p$ of partitions

\item[$ord\phantom{sp}$] : $\ord\mathbf m$

\item[$idx\phantom{sp}$] : $\idx\mathbf m$

\item[$fuchs\!$] : Fuchs relation (no meaning if $\mathbf m$ is a tuple of partitions)

\item[$rod\phantom{sp}$] : the difference of the value $\ord\mathbf m$ under one-step reduction

\item[$redsp$] : the positions $(j_0^{(N)},\dots,j_{p-1}^{(N)})$ in \eqref{eq:generate} 
specified by one-step reduction
\item[$fspt\phantom{p}$] : corresponding fundamental tuple
\end{itemize}
In place of the tuple $\mathbf m$, the corresponding element 
$\alpha\in Q_+$ in the Kac-Moody root space may be indicated by
\begin{equation}\label{eq:nkac}
  \mathbf m=\texttt{[$n$,[$n_{0,1}$,$\ldots$,$n_{0,n_0-1}$],$\ldots$,[$n_{p-1,1}$,$\ldots$,
$n_{p-1,{n_{p-1}-1}}$]]}\quad(\text{cf.~\eqref{eq:n2m}}).
\end{equation}
\begin{itemize}
\item
$\verb|opt="sp"|$ : return spectral type

\item
$\verb|opt="basic"|$ : return $fspt$ (return 0 if $\mathbf m\not\in\Delta_+$)

\item
$\verb|opt="construct"|$ : return the construction from a fundamental tuple
(return 0 if $\mathbf m\not\in\Delta_+$ ({\sl not realizable}))

\item
$\verb|opt="strip"|$ : delete the exponents 0  from the tuple

\item
$\verb|opt="sort"|$ : return the tuple with sorted partitions

\item
$\verb|opt="root"|$ : return the construction by Kac-Moody Weyl group

The returned list has 3 elements which are base spectral type, given spectral type and
the list of reflections. The reflection has 3 elements, namely, the value of inner product,
the branch and the position in the star-shaped Dynkin diagram which specified the root in
the Dynkin diagram.

\item
$\verb|opt="idx"|$ : return index of rigidity

\item
$\verb|opt="kac"|$ : return the corresponding element in $Q_+$ (cf.~\eqref{eq:nkac}).
\end{itemize}

\vspace{5mm}
\hfill\scalebox{0.85}{\begin{tikzpicture}
 [root/.style={draw,circle,inner sep=0mm,minimum size=2.6mm}]
\node[root] (a1) at (0,0) [label=below:$\alpha_{0,2}$, label=above:$1$]{};
\node[root] (a2) at (1,0)  [label=below:$\alpha_{0,1}$, label=above:$3$]{};
\node[root] (a3) at (2,0)  [label=below:$\alpha_0$,label=above:$4\ \ \ \ \ $]{};
\node[root] (a4) at (3,0)  [label=below:$\alpha_{2,1}$,label=above:$3$]{};
\node[root] (a5) at (4,0)  [label=below:$\alpha_{2,2}$,label=above:$2$]{};
\node[root] (a6) at (5,0)  [label=below:$\alpha_{2,3}$,label=above:$1$]{};
\node[root] (a7) at (2,1)  [label=right:$\alpha_{1,1}$,label=left:$2$]{};
\draw (a1)--(a2)--(a3)--(a4)--(a5)--(a6) (a3)--(a7);
\end{tikzpicture}}\\[-26mm]
\noindent
\begin{verbatim}
[0] os_md.chkspt([[1,2,1],[2,2],[1,1,1,1]]);
 [3,4,2,0,1,[ 1 0 0 ],[[1],[1],[1]]]]
[1] os_md.chkspt("121,22,1111");
 [3,4,2,0,1,[ 1 0 0 ],[[1],[1],[1]]]]
[2] os_md.chkspt("121,22,1^4");
 [3,4,2,0,1,[ 1 0 0 ],[[1],[1],[1]]]]
[3] os_md.chkspt([[1,2,1],[2,2],[1,1,1,1]]|opt="kac");
 [4,[3,1],[2],[3,2,1]]
[4] os_md.chkspt([4,[3,1],[2],[3,2,1]]);
 [3,4,2,0,1,[ 1 0 0 ],[[1],[1],[1]]]]
[5] os_md.chkspt([4,[3,1],[2],[3,2,1]]|opt="sp");
 [[1,2,1],[2,2],[1,1,1,1]]
[6] os_md.chkspt([[1,2,1],[2,2],[1,1,1,1]]|opt="sort");
 [[2,1,1],[2,2],[1,1,1,1]]
[7] os_md.chkspt([[1,2,1],[2,2],[1,1,1,1]]|opt="basic");
 [[1],[1],[1]]
[8] os_md.chkspt([[1,2,1],[2,2],[1,1,1,1]]|opt="construct");
 [
  [[1],[1],[1]],
  [[1,1],[1,1],[1,1]],
  [[1,1,1],[2,1],[1,1,1]],
  [[2,1,1],[2,2],[1,1,1,1]]
 ]
[9] os_md.chkspt([[0,1,2,1],[2,0,2,0],[1,1,1,1]]|opt="strip");
 [[2,1,1],[2,2],[1,1,1,1]]
[10] os_md.chkspt("21,21,21,21"|opt="root");
 [[[1],[1],[1],[1]],[[2,1],[2,1],[2,1],[2,1]],
 [[1,3,1],[1,2,1],[1,1,1],[1,0,1],[2,0,0]]]
[11] os_md.chkspt("21,21,21,111");
 [4,3,0,0,1,[ 0 0 0 0 ],[[1,1],[1,1],[1,1],[1,1]]]
[12] os_md.chkspt([[3,2],[2,2,1],[1,1,1,1]]);
 illegal partitions
 -1
[13] os_md.chkspt("31,31,31,22");
 not realizable
 0
[14] os_md.chkspt("43,322,1^7");
 [3,7,0,0,1,[ 0 0 0 ],[[3,3],[2,2,2],[1,1,1,1,1,1]]]
[15] os_md.chkspt("21,21,21,111");
 [4,3,0,0,1,[ 0 0 0 0 ],[[1,1],[1,1],[1,1],[1,1]]]
[16] os_md.chkspt([[2*m,2*m],[m,m,m,m],[m,m,m,m-1,1]]|opt="idx");
 -2*m+2
\end{verbatim}

\noindent
The meaning of the above results is as follows.

\texttt{[0]}--\texttt{[2]} : tuple $\mathbf m=121,22,1111$ is checked.

$\Rightarrow$ the result \texttt{[3,4,2,0,1,[ 1 0 0 ],[[1],[1],[1]]]]} means

\quad $\# $ partitions : 3

\quad $\ord \mathbf m=4$

\quad $\idx\mathbf m=2$

\quad Fuchs condition : no meaning

\quad The order decreased by one-step reduction : 1

\quad The position for the reduction : $(\nu_0,\nu_1,\nu_2) = (1,0,0)$ in \eqref{eq:red1}.  

\quad It is in the orbit of \ "1,1,1"$=$\texttt{[[1],[1],[1]]}$=\alpha_0$.

\texttt{[3]} : the expression \eqref{eq:nkac}  is returned:

\quad \texttt{[[$n$],[$n_{0,1}.n_{0,2},\ldots$],[$n_{1,1},\ldots$],$\ldots$]} 

\texttt{[8]} : the construction \eqref{eq:generate} is given.

\texttt{[10]} : the construction of $21,21,21,21=3\alpha_0+\sum\limits_{j=0}^3 \alpha_{j,1}$ 
(cf.~\eqref{eq:generate}) is given:

\vspace{-0.3mm}
\quad It is constructed from $1,1,1,1=\alpha_0$ by simple refrections $s_\alpha$ (cf.~\eqref{eq:sref}).

\quad $[1,3,1]$ means that $s_{\alpha_{3,1}}$ increase the coefficient of $\alpha_{3,1}$ 
by $1$.

\texttt{[13]} : the answer means that $31,31,32,22$ does not correspond to a root.

\medskip
We can consider a similar problem in the case of a general 
symmetrizable Kac-Moody root system
which includes the star-shaped Kac-Moody root system. 
The imaginary fundamental root $\alpha$ is characterized by
$(\alpha|\gamma)\le 0$ for $\gamma\in\Pi$, where $\Pi$ is the set of simple roots
defining the root system.
It is proved by \cite{HO} that
for a given integer $N$, there exist only finite 
fundamental roots $\alpha$ with a given index of rigidity $(\alpha|\alpha)=N$ 
under a simple equivalence relation and classified
them when $N\ge -2$.  Note that the imaginary fundamental root is called a 
basic root in \cite{HO}.
It is desirable to have a computer program to get the fundamental roots $\alpha$
with a given index of rigidity. 

The fundamental roots of the symmetrizable Kac-Moody root system
corresponds to the spectral types of ordinary differential equations 
\begin{equation*}
 \mathcal N : \frac {du}{dx}=\Bigl(\sum_{j=1}^{q-1}\sum_{k=1}^{m_k}\frac{A_{j,k}}{(x-a_j)^k}
 +\sum_{k=2}^{m_q}A_{q,k}x^{k-2}\Bigr)u\qquad(A_{j,k}\in M(n,\mathbb C))
\end{equation*}
which allow unramified 
irregular singularities together with regular singularities (cf.~\cite{Over}).
The number of accessory parameters of the equation equals $2-(\alpha|\alpha)$.
The deformation preserving the monodromy of solutions of the equation gives 
$\bigl(2-(\alpha|\alpha)\bigr)$-dimensional Painlev\'e-type equation and 
it is known that the Painlev\'e-type equation is determined by $\widetilde W\alpha$.
When $(\alpha|\alpha)=0$, it corresponds to a classical Painlev\'e equation. 
Hence the classification of fundamental roots are closely related to the classification of
higher dimensional Painlev\'e-type equations.

The classification of fundamental roots $\alpha$ with $(\alpha|\alpha)=-2$ 
in \cite{HO} induces 
the classification of 4-dimensional Painlev\'e-type equations in \cite{HS} and 
that with $(\alpha|\alpha)=-4$ in \cite{Orims,Ow} induces 
8 new 6-dimensional Painlev\'e-type equations constructed by \cite{Su}.

Note that the differential equation $\dot u=A(x)u$ with a matrix-valued rational function $A(x)$ 
can be written as above.
As in the case of Fuchsian ordinary differential equations, the spectral type of the equation
with unramified irregular singularities can be also defined through the versal unfolding of the 
singularities into regular singularities (cf.~\cite{Over}).

The spectral type of the unfolded Fuchsian equation $\mathcal M$ is expressed by 
a tuple of partitions corresponding to a root of star-shaped Kac-Moody root system and 
that of the original equation $\mathcal N$ is a refinement of partitions in the tuple.
Then $p=m_1+\cdots+m_q$ and $m_j\ge 1$.
The spectral types of such equation and the refinements of realizable tuples of
the partitions are in one to one correspondence, which is conjectured by the author 
\cite{Okzv, Over} and recently proved by K.~Hiroe \cite{Hi}.
Hence the spectral type of the equation $\mathcal N$ with  unramified irregular singularities 
is obtaineed by a refinement of a Fuchsian differential equation corresponding to  
the confluence of regular singularities. 

For two partitions of 6 given by $6=4+2=2+2+1+1$, there are refinements 
$6=4(=2+2)+2(=1+1)$ and $6=4(=2+1+1)+2(=2)$ and describe 
spectral types of confluences of two regular singular points 
to an irregular singular point.
They are expressed by 
$2211|42$ and $2112|42$, respectively. 
When some partitions in a spectral type 
form a refinement, they are connected with  $|$ 
in the order according to refinement and then the order of integers in this 
expression is essential.
They may be simply written as $(22)(11)$ and $(211)(2)$, respectively (cf.~\cite{HO}).

Middle convolutions and additions are similarly defined and they induce similar 
transformations of tuples of partitions \eqref{def:sinu} and \eqref{def:s0} which are 
compatible to the structure of their refinements.

The spectral types 
\[11|11,11,11\quad 11|11|11,11\quad 11|11,11|11\quad 11|11|11|11\]
are confluences of the spectral type $11,11,11,11$ of the 2-nd order Fuchsian equation 
$\mathcal M$ $(p=4,\ n=2)$ with 4 singular points.  
The deformations of the corresponding equations $\mathcal N$ $(m_j\ge 1,\ m_1+\cdots+m_q=4,\ n=2)$
with these spectral types and $11,11,11,11$ which preserve the monodromies 
give Painlev\'e equations of type V, IV, III, II and VI,
respectively.

\medskip
The function \texttt{refinements()} in 
\cite{ORisa} gives all the refinements of a given tuple.

\smallskip
\noindent
\texttt{refinements}%
\texttt{($l$}\verb/|/\texttt{opt=$f$)}\\
:: get all the refinements for several partitions of a positive integer

\begin{itemize}
\item
\texttt{opt="s"} : Get all the unramified irregular spectra by the confluence of a given Fuchsian spectrum 
by strings using ``\verb/, |/" (cf.~\cite{Over})
\item
\texttt{opt="()"} :Get all the unramified irregular spectra by the confluence of a given Fuchsian spectrum 
by strings using ``\verb/, ( )/" (cf.~\cite{HS})
\end{itemize}
\begin{verbatim}
 [0] R=os_md.refinements("11,11,11,11"|opt="s")$
 [1] for(T=R;T!=[];T=cdr(T)) print(car(T));
  11,11,11,11
  11|11|11|11
  11,11|11|11
  11,11,11|11
  11|11,11|11
 [2] R=os_md.refinements("11,11,11,11"|opt="()")$
 [3] for(T=R;T!=[];T=cdr(T)) print(car(T));
  1 1,1 1,1 1,1 1
  (((1))) (((1)))
  1 1,((1)) ((1))
  1 1,1 1,(1) (1)
  (1) (1),(1) (1)
 [4] R=os_md.refinements("42,21111,21111"|opt="()")$
 [5] os_md.mycat(R|delim="\n");
  2 1 1 1 1,2 1 1 1 1,4 2
  ((1) (1) (1) (1)) ((2))
  ((2) (1) (1)) ((1) (1))
  4 2,(2) (1) (1) (1) (1)
  2 1 1 1 1,(1 1 1 1) (2)
  2 1 1 1 1,(2 1 1) (1 1)
\end{verbatim}
\medskip
Lastly we show a function analyzing equation with a given spectral type.

\medskip
\noindent
\texttt{chkcspt}%
$\verb|(|m\verb/|show=1)/$\\
:: analyze spectral type of ODE with unramified irregular singularities
\begin{itemize}
\item
spectral type $m$ is indicated by a list or a string
\item
Return the list 
\texttt{[}number of singular points, Poincare ranks at singular points, rank of $m$, index of rigidity, 
maximal rank decreased by a reduction, places where the multiplicities are decreased by the 
reduction, original spectral type, reduced spectral type\texttt{]}

See an example with the option \texttt{show=1}.
\end{itemize}
\begin{verbatim}
 [0] os_md.chkcspt("1111|211,22");
  [2,[1,0],4,2,1,[0,0],
  1111|211,22, [[[2,[1,1]],[1,[1]],[1,[1]]],[2,2]],
  111|111,12, [[[1,[1]],[1,[1]],[1,[1]]],[1,2]]]
 [1] os_md.chkcspt("1111|211,22"|show=1);
  1111|211,22   (1 1) (1) (1),2 2
  points:    2  with Poincare ranks  [1,0]
  rank:      4
  index:	    2
  reduct:    1 at [0,0] -> 111|111,12   (1) (1) (1),1 2
 [2] os_md.chkcspt("1111|211|22"|show=1);
  1111|211|22   ((1 1)) ((1) (1))
  points:    1  with Poincare ranks  [2]
  rank:      4
  index:     2
  reduct:    1 at [0] -> 111|111|12   ((1)) ((1) (1))
\end{verbatim}

\bigskip
\centerline{\bf Acknowledgements}

\smallskip
In 2008, a program \texttt{okubo} \cite{okubo} was written by the author using the computer 
language \texttt{C}.
It has several functions handling the tuples of partitions related to the spectral types of
Fuchsian differential equations including 
a function which generates fundamental tuples with a given index.

In 2019, H.~Kawakami pointed out a fundamental tuple with the index $-4$ which 
is not generated by
the function (cf.~Remark~\ref{rem:fund} (ii)) and then the author 
rewrote the function \texttt{spbasic(\,)} using the computer algebra 
\texttt{Risa/Asir}.

In 2024, L.~Yongchao showed the author a missing fundamental tuple
with the index $-18$ and pointed out 
a connection to a problem in mathematical physics (cf.~\cite{Mphy}).
It is the only one missing fundamental tuple when the index of rigidity $\ge-24$ and there exists another one in 2889 fundamental tuples with the index $= -26$.
The author has rewritten \texttt{spbasic(\,)} and writes this paper explaining the
algorithm used by \texttt{spbasic(\,)}.

The author greatly thanks to H.~Kawakami and L.~Yongchao for pointing out missing fundamental 
tuples caused by the original programs.


\begin{thebibliography}{GKZ}
\bibitem[CB]{CB}
Crawley-Boevey W.,
On matrices in prescribed conjugacy classes with no common invariant subspaces,
and sum zero,
Duke Math.~J. \textbf{118}(2003), 339-352.
%
\bibitem[Hi]{Hi}
Hiroe K,
Deformation of moduli spaces of meromorphic G-connections on P1 via unfolding of irregular singularities,
arXiv.2407.20486, 2024.
%
\bibitem[HKNS]{HS}
Hiroe K, H. Kawakami, A. Nakamura and H. Sakai, 
\textsl{4-dimensional Painlev\'e-type equations},
MSJ Memoirs \textbf{37}, 
Mathematical Society of Japan, 2018.

\bibitem[HO]{HO}
Hiroe K.\ and Oshima T.,
{Classification of roots of symmetric Kac-Moody root systems and its application, 
Symmetries}, Integrable Systems and Representations, Springer Proceedings in Mathematics and 
Statistics \textbf{40} (2012), 195--241.
\bibitem[Kc]{Kc}
Kac V.~C.,
\textsl{Infinite dimensional Lie algebras}, Third Edition, Cambridge Univ.\ Press 1990.
\bibitem[Ko]{Ko}
Kostov,~V.~P., 
{The Deligne-Simpson problem for zero index of rigidity},
\textsl{Perspective in Complex Analysis}, 
Differential Geometry and Mathematical Physics, 
{World Scientific} 2001, 1--35.
\bibitem[MOTZ]{Mphy}
Mekareeya1 N., Ohmori K., Yuji Tachikawa Y.\ and Zafrir G.,
{$E_8$ instantons on type-A ALE spaces
and supersymmetric field theories},
J.\ High Energ.\ Phys.\ \textbf{144} (2017).
\bibitem[O1]{Orims}
Oshima T.,
{Classification of Fuchsian systems and their connection problem}, 
RIMS K\^oky\^uroku Bessatsu {\bf B37} (2013), 
163--192.
\bibitem[O2]{Ow}
Oshima T.,
\textit{Fractional calculus of Weyl algebra and Fuchsian differential equations},
MSJ Memoirs \textbf{28}, Mathematical Society of Japan, Tokyo, 2012.
\bibitem[O3]{Ored}
Oshima T.,
{Reducibility of hypergeometric equations}, 
Analytic, Algebraic and Geometric Aspects of
Differential Equations, Trends in Mathematics, Springer,
429--453, 2017.
\bibitem[O4]{Okzv}
Oshima, T.,
{Confluence and versal unfolding of {P}faffian systems}, 
Josai Mathematical Monographs \textbf{12} (2020), 117--151.

\bibitem[O5]{Over}
Oshima, T.,
{Versal unfolding of irregular singularities of a linear differential equation on 
the Riemann sphere},
Publ. RIMS Kyoto Univ.~\textbf{57}, 
(2021). 893--920.

\bibitem[O6]{Ofund}
Oshima T.,
List of fundamental spectral types with index of the rigidity $\ge-10$,
2024,
\texttt{https://www.ms.u-tokyo.ac.jp/\~{}oshima/paper/spect10.pdf}

\bibitem[O7]{okubo}
Oshima T.,
\texttt{okubo},
a computer program for Katz/Yokoya/Oshima algorithm, 2007--2008,
\texttt{http://www.ms.u-tokyo.ac.jp/\~{}oshima/}
\bibitem[O8]{ORisa}
Oshima T.,
\texttt{os\_muldif.rr}, a library of computer algebra \texttt{Risa/Asir},
2008--2024.\\
\url{http://www.ms.u-tokyo.ac.jp/\~oshima/}

\bibitem[Risa]{Risa}
Noro M.\ etc.,
\texttt{Risa/Asir}, Computer algebra, \\
\texttt{http://www.math.kobe-u.ac.jp/Asir/asir.html}
\url{http://www.math.kobe-u.ac.jp/Asir/asir.html}

\bibitem[Su]{Su}
Suzuki T.,
Six-dimensional Painlev\'e systems and their particular solutions in terms of rigid systems,
J.\ Math.\ Physics \texttt{55} (2014), 57--69.
\end{thebibliography}
\end{document}